\documentclass[11pt,a4paper,leqno,twoside]{amsart}
\usepackage{latexsym,amssymb,amsmath}
\input xy
\xyoption{all}

\def\trace{{\rm trace}}

\def\cE{{\mathcal E}}

\def\cN{{\mathcal N}}

\def\cP{{\mathcal P}}

\def\cO{{\mathcal O}}

\def\CC{\mathbb C}
\def\SP{{\rm Spin}_{10}}
\def\RR{\mathbb R}
\def\HH{\mathbb H}
\def\AA{{\mathbb A}}

\def\OO{\mathbb O}

\def\SS{\mathbb S}

\def\11{\mathbf 1}
\def\PP{\mathbb P}
\def\QQ{\mathbb Q}

\def\e1{\varepsilon_1}
\def\e2{\varepsilon_2}
\def\e3{\varepsilon_3}

\def\P2{{\PP}^2}

\def\JO{{\mathcal J}_3(\OO)}
\def\JZ{{\mathcal J}_2(\OO)}
\def\00{\underline{0}}
\def\J0{{\cal J}_3(\underline{0})}

\def\PJ0{\PP({\cal J}_3(\underline{0}))}

\def\fe{{\mathfrak e}}

\def\a{\alpha}

\def\om{\omega}

\def\g{\gamma}
\def\s{\sigma}

\def\e{\varepsilon}

\def\op{{\mathord{\,\oplus }\,}}

\def\ra{{\mathord{\;\rightarrow\;}}}

\def\AP2{{\AA\PP}^2}
\def\RP2{{\RR\PP}^2}
\def\CP2{{\CC\PP}^2}
\def\HP2{{\HH\PP}^2}
\def\OP2{{\OO\PP}^2}

\def\dim{{\rm dim}\;}

\newtheorem{theo}{Theorem}[section]

\newtheorem{lemm}[theo]{Lemma}
\newtheorem{prop}[theo]{Proposition}

\begin{document}

\title{The Chow ring of the Cayley plane}
\author{A. Iliev, L. Manivel}

\begin{abstract} We give a full description of the Chow ring of the complex Cayley 
plane $\OO\PP^2$. For this, we describe explicitely the most interesting of its Schubert 
varieties and compute their intersection products. Translating our results in the 
Borel presentation, i.e. in terms of Weyl group invariants, we are able to compute
the degree of the variety of reductions $Y_8$ introduced in \cite{IM}.
\end{abstract}

\maketitle

\section{Introduction}

In this paper we give a detailed description of the Chow ring of the complex 
Cayley plane\footnote{Not to be confused with the real Cayley plane $F_4/Spin_9$,
the real part of $\OO\PP^2$, which admits a cell decomposition $\RR^0\cup \RR^8\cup \RR^{16}$
and is topologically much simpler.}
 $X_8=\OO\PP^2$, the fourth Severi variety. This is a smooth complex projective
variety of dimension $16$, homogeneous under the action of the adjoint group of 
type $E_6$. It can be described as the closed orbit in the projectivization $\PP^{26}$ 
of the minimal representation of $E_6$. 

The Chow ring of a projective homogeneous variety $G/P$ has been described classically
in two different ways. 

First, it can be described has a quotient of a ring of invariants. Namely, we have to 
consider the action of the Weyl group of $P$ on the character ring, take the invariant
subring, and mod out by the homogeneous ideal generated by the invariants (of positive 
degree) of the full Weyl group of $G$. This is the {\sl Borel presentation}. 

Second, the Chow ring has a basis given by the {\it Schubert classes}, the classes of the 
closures of the $B$-orbits for some Borel subgroup $B$ of $E_6$. These varieties
are the Schubert varieties. Their intersection products can in principle be 
computed by using Demazure operators \cite{bgg}. This is  the {\sl Schubert presentation}. 

We  give a detailed description of the Schubert presentation of the Chow ring 
$A^*(\OO\PP^2)$ of the Cayley plane. We describe explicitely the most interesting Schubert 
cycles, after having explained how to understand geometrically a Borel subgroup of $E_6$. 
Then we compute the intersection numbers. In the final section, we turn to the Borel presentation 
and determine the classes of some invariants of the partial Weyl group in terms of Schubert
classes, from which we deduce the Chern classes of the normal  bundle of $X_8=\OO\PP^2$ in 
$\PP^{26}$. This allows us to compute the degree of the variety of reductions $Y_8\subset 
\PP^{272}$ introduced in \cite{IM}, which was the initial motivation for writing 
this note.

\section{The Cayley plane}

Let ${\bf O}$ denote the normed algebra of (real) octonions, and let 
$\OO$ be its complexification. The space
$$\JO = \Bigg\{ 
\begin{pmatrix} 
c_1            & x_3            & \bar{x}_2  \\ 
\bar{x}_3 & c_2            &  x_1            \\ 
x_2            & \bar{x}_1 &  c_3            \\ 
\end{pmatrix}
: c_i \in {\bf C}, x_i \in \OO \Bigg\} \cong {\bf C}^{27}
$$
of $\OO$-Hermitian matrices of order $3$, is the exceptional 
simple complex Jordan algebra, for the Jordan multiplication 
$A \circ B = \frac{1}{2}(AB+BA)$. 

The subgroup $SL_3(\OO)$ of $GL(\JO)$ consisting in automorphisms preserving the 
determinant is the adjoint group of type $E_6$. The Jordan algebra $\JO$ and its
dual are the minimal representations of this group. 

The action of $E_6$ on the 
projectivization $\PP\JO$ has exactly three orbits: the complement of the 
determinantal hypersurface, the regular part of this hypersurface, and its singular
part which is the closed $E_6$-orbit. These three orbits are the sets of matrices 
of rank three, two, and one respectively. 

The closed orbit, i.e. the (projectivization of) the set of rank one matrices,
is the {\it Cayley plane}. It can be defined by the quadratic equation 
$$X^2=\trace(X)X, \qquad  X\in \JO,$$ 
or as the closure of the affine cell 
$$\OO\PP^2_1=
\Bigg\{\begin{pmatrix} 1 & x & y \\ \bar{x} & x\bar{x} & y\bar{x} \\ 
\bar{y} & x\bar{y} & y\bar{y} \end{pmatrix}, \quad x,y\in\OO\Bigg\}\cong {\bf C}^{16}.$$
It is also the closure of the two similar cells
$$\OO\PP^2_2=
\Bigg\{\begin{pmatrix} \bar{u}u & u & vu \\ \bar{u} & 1 & v \\ 
\bar{u}\bar{v} & \bar{v} & v\bar{v} \end{pmatrix}, \quad u,v\in\OO\Bigg\}\cong {\bf C}^{16},$$
$$\OO\PP^2_3=
\Bigg\{\begin{pmatrix} \bar{t}t & \bar{s}t & t \\ \bar{t}s & \bar{s}s & s \\ 
\bar{t} & \bar{s} & 1 \end{pmatrix}, \quad s,t\in\OO\Bigg\}\cong {\bf C}^{16}.$$
Unlike the ordinary projective plane, these three affine cells do not cover $\OO\PP^2$. 
The complement of their union is 
$$\OO\PP^2_{\infty}=
\Bigg\{\begin{pmatrix} 0 & x_3 & x_2 \\ \bar{x_3} & 0 & x_1 \\ 
\bar{x_2} & \bar{x_1} & 0 \end{pmatrix}, \quad \begin{matrix} q(x_1) = q(x_2) = 
q(x_3) = 0, \\ 
x_2\bar{x_3} = x_1x_3 = \bar{x_1}x_2 = 0\end{matrix}   \Bigg\},$$
a singular codimension three linear section. Here, $q(x)=x\bar{x}$ denotes the non degenerate 
quadratic form on $\OO$ obtained by complexification of the norm of ${\bf O}$. 

\medskip
Since the Cayley plane is a closed orbit of $E_6$, it can also be identified with the 
quotient of $E_6$ by a parabolic subgroup, namely the maximal parabolic subgroup defined 
by the simple root $\a_6$ in the notation below. The semi-simple part 
of this maximal parabolic is isomorphic to ${\rm Spin}_{10}$.
\medskip

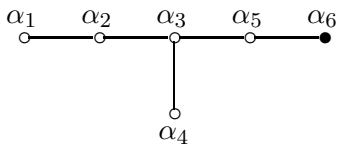
\begin{center}
\setlength{\unitlength}{5mm}
\begin{picture}(20,5)(-5,0)
\multiput(-.3,3.8)(2,0){4}{$\circ$}
\multiput(0,4)(2,0){4}{\line(1,0){1.7}}
\put(3.7,1.8){$\circ$}
\put(7.7,3.8){$\bullet$}
\put(-.6,4.4){$\a_1$}
\put(1.4,4.4){$\a_2$}
\put(3.4,4.4){$\a_3$}
\put(5.4,4.4){$\a_5$}
\put(7.4,4.4){$\a_6$}
\put(3.45,1.3){$\a_4$}
\put(3.85,2.1){\line(0,1){1.8}}
\end{picture}
\end{center}

\medskip
The $E_6$-module $\JO$ is {\it minuscule}, meaning that its weights with respect 
to any maximal torus of $E_6$, are all conjugate under the Weyl group action. 
We can easily list these weights as follows. Once we have fixed a set of simple 
roots of the Lie algebra, we can define the height of any weight $\omega$ as the sum of 
its coeficients when we express $\omega$ on the basis of simple roots. Alternatively, 
this is just the scalar product $(\rho ,\omega)$, if $\rho$ denotes, as usual, the sum 
of the fundamental weights, and the scalar product is dual to the Killing form. 
The highest weight $\om_6$ of $\JO$ is the unique weight with maximal height. 
We can obtain the other weights using the following process: if we have some 
weight $\om$ of $\JO$, we express it in the basis of fundamental weights. For each 
fundamental weight $\om_i$ on which the coefficient of $\om$ is positive, we apply the
corresponding simple reflection $s_i$. The result is a weight of $\JO$ of height
smaller than that of $\om$, and we obtain all the weights in this way. The following 
diagram is the result of this process. We do not write down the weights explicitely, 
but we keep track of the action of the simple reflections: if we apply $s_i$ to go from 
a weight to another one, we draw an edge between them, labeled  with an $i$. 

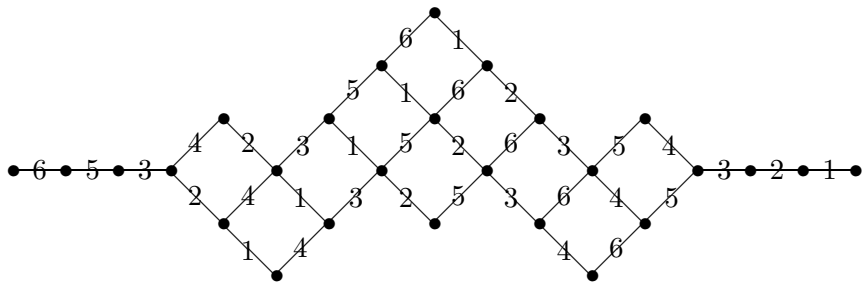
\medskip

\begin{center}
\setlength{\unitlength}{3.5mm}
\begin{picture}(40,12)(-8,2)

\put(-6,8){\line(1,0){6}}
\put(0,8){\line(1,1){2}}
\put(0,8){\line(1,-1){4}}

\put(2,10){\line(1,-1){4}}
\put(2,6){\line(1,1){8}}
\put(4,4){\line(1,1){8}}
\put(8,12){\line(1,-1){8}}
\put(10,14){\line(1,-1){8}}
\put(6,10){\line(1,-1){4}}
\put(10,6){\line(1,1){4}}
\put(14,6){\line(1,1){4}}
\put(16,4){\line(1,1){4}}
\put(18,10){\line(1,-1){2}}
\put(20,8){\line(1,0){6}}
\multiput(-6.2,7.7)(2,0){4}{$\bullet$}
\multiput(1.8,9.7)(4,0){5}{$\bullet$}
\multiput(1.8,5.7)(4,0){5}{$\bullet$}
\multiput(3.8,7.7)(4,0){5}{$\bullet$}
\multiput(3.8,3.7)(12,0){2}{$\bullet$}
\multiput(7.8,11.7)(4,0){2}{$\bullet$}
\multiput(21.8,7.7)(2,0){3}{$\bullet$}
\put(9.8,13.7){$\bullet$}
\put(-5.2,7.65){$6$}
\put(-3.2,7.65){$5$}
\put(-1.2,7.65){$3$}
\multiput(.7,6.7)(2,2){2}{$2$}
\multiput(.7,8.7)(2,-2){3}{$4$}

\multiput(2.7,4.6)(2,2){5}{$1$}
\multiput(4.8,8.6)(2,-2){2}{$3$}
\multiput(8.7,6.6)(2,2){3}{$2$}
\multiput(6.7,10.7)(2,-2){3}{$5$}
\multiput(8.7,12.7)(2,-2){5}{$6$}

\multiput(12.7,6.6)(2,2){2}{$3$}
\multiput(14.7,4.6)(2,2){3}{$4$}
\multiput(16.8,8.6)(2,-2){2}{$5$}
\put(20.8,7.65){$3$}\put(22.8,7.65){$2$}\put(24.8,7.65){$1$}
\end{picture}
\end{center}

\section{The Hasse diagram of Schubert cycles}

Schubert cycles in $\OO\PP^2$ are indexed by a subset $W^0$ of the Weyl group
$W$ of $E_6$, the elements of which are minimal length representatives of the 
$W_0$-cosets in $W$. Here $W_0$ denotes the Weyl group of the maximal parabolic 
$P_6\in E_6$: it is the subgroup of $W$ generated by the simple reflections 
$s_1,\ldots ,s_5$, thus isomorphic to the Weyl group of $\SP$.

But $W_0$ is also the stabilizer in $W$ of the weight $\om_6$. Therefore, the
weights of $\JO$ are in natural correspondance with the elements of $W^0$, 
and we can obtain very explicitely,  from the picture above,  
the elements of $W^0$. Indeed, choose
any vertex of the diagram, and any chain of minimal length joining this vertex 
to the leftmost one. Let $i_1,\ldots ,i_k$ be the consecutive labels on the edges
of this chain; then $s_{i_1}\cdots s_{i_k}$ is a minimal decomposition of the 
corresponding elements of $W^0$, and every such decomposition is obtained in this way. 
 \smallskip

For any $w\in W^0$, denote by $\s_w$ the corresponding Schubert cycle of $\OO\PP^2$. 
This cycle  $\s_w$  belongs to $A^{l(w)}(\OP2)$, where $l(w)$ denotes the length of $w$. 
We have just seen that this length is equal to the distance of the point corresponding 
to $w$ in the picture above, to the leftmost vertex. In particular, the dimension 
of $A^k(\OP2)$ is equal to one for $0\le k\le 3$, to two for $4\le k\le 7$, to three 
for $k=8$ (and by duality, this dimension is of course unaltered when $k$ is changed 
into $16-k$). 

The degree of each Schubert class can be deduced from the Pieri formula, which is
particularly simple in the minuscule case. Indeed, we have (\cite{hiller}, Corollary
3.3), if $H$ denotes the hyperplane class:
$$\s_w.H^k=\sum_{l(v)=l(w)+k}\kappa (w,v)\s_v,$$
where $\kappa(w,v)$ denotes the number of path from $w$ to $v$ in the diagram
above; that is, the number of chains $w=u_0\ra u_1\ra\cdots\ra u_k=v$ in $W^0$
such that $l(u_i)=l(w)+i$ and $u_{i+1}u_i^{-1}$ is a simple reflection. In particular, 
the degree of $\s_w$ is just $\kappa(w,w^0)$, where $w^0$ denotes the longest element 
of $W^0$, which corresponds to the leftmost vertex of the diagram. We include these
degrees in the following picture, the Hasse diagram of $\OO\PP^2$. Note that they can
very quickly be computed inductively, beginning from the left: the degree of each cycle
is the sum of the degrees of the cycles connected to it in one dimension less.

\begin{center}
\setlength{\unitlength}{3.5mm}
\begin{picture}(40,12)(-8,2)

\put(-6,8){\line(1,0){6}}
\put(0,8){\line(1,1){2}}
\put(0,8){\line(1,-1){4}}
\put(2,10){\line(1,-1){4}}
\put(2,6){\line(1,1){8}}
\put(4,4){\line(1,1){8}}
\put(8,12){\line(1,-1){8}}
\put(10,14){\line(1,-1){8}}
\put(6,10){\line(1,-1){4}}
\put(10,6){\line(1,1){4}}
\put(14,6){\line(1,1){4}}
\put(16,4){\line(1,1){4}}
\put(18,10){\line(1,-1){2}}
\put(20,8){\line(1,0){6}}
\multiput(-6.2,7.7)(2,0){4}{$\bullet$}
\multiput(1.8,9.7)(4,0){5}{$\bullet$}
\multiput(1.8,5.7)(4,0){5}{$\bullet$}
\multiput(3.8,7.7)(4,0){5}{$\bullet$}
\multiput(3.8,3.7)(12,0){2}{$\bullet$}
\multiput(7.8,11.7)(4,0){2}{$\bullet$}
\multiput(21.8,7.7)(2,0){3}{$\bullet$}
\put(9.8,13.7){$\bullet$}
\multiput(-6.2,6.7)(2,0){4}{$1$}
\put(1.8,4.7){$1$}
\put(1.8,10.4){$1$}
\put(3.8,2.7){$1$}
\put(3.8,8.5){$2$}
\put(5.8,4.7){$3$}
\put(5.8,10.5){$2$}
\put(7.8,6.7){$5$}
\put(7.8,12.5){$2$}
\put(9.8,4.7){$5$}
\put(9.8,8.7){$7$}
\put(9.8,14.5){$2$}
\put(11.5,6.7){$12$}
\put(11.8,12.5){$9$}
\put(13.5,4.7){$12$}
\put(13.8,10.5){$21$}
\put(15.5,2.7){$12$}
\put(15.5,8.6){$33$}
\multiput(19.5,6.7)(2,0){4}{$78$}
\put(17.5,4.7){$45$}
\put(17.5,10.4){$33$}
\end{picture}
\end{center}

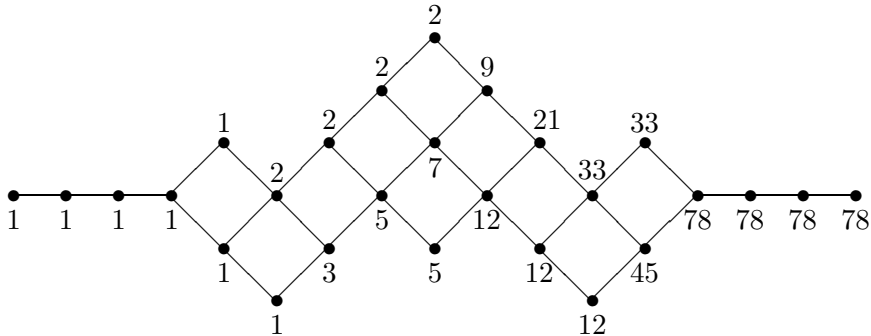
\centerline{{\it Degrees of the Schubert cycles}}

\bigskip
We can already read several interesting informations on this diagram. 
\begin{enumerate}
\item The degree of $\OO\PP^2\subset\PP^{26}$ is $78$. This is precisely the 
dimension of $E_6$. Is there a natural explanation of this coincidence ?  
\item One of the three Schubert varieties of dimension $8$ is a quadric. 
This must be an $\OO$-line in $\OO\PP^2$, i.e. a copy of $\OO\PP^1\simeq\QQ^8$. 
Indeed, $E_6$ acts transitively on the family of  these lines, which is actually 
parametrized by $\OO\PP^2$ itself. In particular, a Borel subgroup has a fixed 
point in this family, which must be a Schubert variety.  
\item The Cayley plane contains two families of Schubert cycles which are 
maximal linear subspaces: a family of $\PP^4$'s, which are maximal linear subspaces
in some $\OO$-line, and a family of $\PP^5$'s which are not contained in any 
$\OO$-line. We thus recover the results of \cite{LMhom}, from which we also know that 
these two families of linear spaces in $\OO\PP^2$ are homogeneous. Explicitely, we
can describe both types in the following way. 

Let $z\in \OO$ be a non zero octonion such that $q(z)=0$. Denote by $R(z)$ and $L(z)$ 
the spaces of elements of $\OO$ defined as the images of the right and left multiplication 
by $z$, respectively. Similarly, if $l\subset\OO$ is an isotropic line,  
denote by $R(l)$ and $L(l)$ the spaces $R(z)$ and $L(z)$, if $z$ is a generator of $l$. 
When $l$ varies, $R(l)$ and $L(l)$ describe the two families of maximal isotropic subspaces 
of $\OO$ (this is a geometric version of triality, see e.g. \cite{chaput3}). 
Consider the sets
$$\begin{array}{r}
\Bigg\{\begin{pmatrix} 1 & x & y \\ \bar{x} & 0 & 0 \\ \bar{y} & 0 & 0 \end{pmatrix},
\; y\in l, x\in L(l)\Bigg\} \quad {\rm and} \quad   
\Bigg\{\begin{pmatrix} 1 & x & 0 \\ \bar{x} & 0 & 0 \\ 0 & 0 & 0 \end{pmatrix},
\; x\in R(l)\Bigg\}\end{array}.$$
Their closures in $\PP\JO$ are maximal linear subspaces of $\OO\PP^2$ of respective 
dimensions $5$ and $4$. \bigskip

\end{enumerate}

\section{What is a Borel subgroup of $E_6$ ?}

The Schubert varieties in $\OO\PP^2$, by definition, are the closures of the $B$-orbits, 
where $B$ denotes a Borel subgroup of $E_6$. To identify the Schubert varieties 
geometrically, we need to understand these Borel subgroups better. 

The Cayley plane $\OO\PP^2= E_6/P_6\subset\PP\JO$ is one of the $E_6$-grassmannians, 
if we mean by this a quotient of $E_6$ by a maximal parabolic subgroup. It is isomorphic
to the dual plane ${\Check \OO\PP}^2=E_6/P_1\subset {\Check \PP\JO}$, the closed 
orbit of the projectivized dual representation. By \cite{LMhom}, we can identify 
$E_6/P_5$ and $E_6/P_3$ to the varieties $G(\PP^1,\OO\PP^2)$ and $G(\PP^2,\OO\PP^2)$ 
of projective lines and planes contained in $\OO\PP^2$. Similarly, $E_6/P_2$ and $E_6/P_3$
can be interpreted as the varieties of projective lines and planes contained in 
${\Check \OO\PP}^2$.

The remaining $E_6$-grassmannian $E_6/P_4$ is the adjoint variety $E_6^{ad}$, the 
closed orbit in the projectivization $\PP\fe_6$ of the adjoint representation. 
By \cite{LMhom} again, $E_6/P_3$ can be identified to the variety $G(\PP^1,E_6^{ad})$ 
of projective lines contained in $E_6^{ad}\subset\PP\fe_6$. \medskip

Now, a Borel subgroup $B$ in $E_6$ is the intersection of the maximal parabolic subgroups
that it contains, and there is one such group for each simple root. Each of these maximal
parabolics can be seen as a point on an $E_6$-grassmannian, and the fact that these 
parabolic subgroups have a Borel subgroup in common, means that these points are {\it incident}
in the sense of Tits geometries \cite{tits}. 

Concretely, a point of $E_6/P_3$ defines a projective plane $\Pi$ in $\OO\PP^2$, a dual
plane ${\Check \Pi}$ in ${\Check \OO\PP}^2$, and a line $\Lambda$ in $E_6^{ad}$. Choose 
a point $p$ and a line $\ell$ in $\OO\PP^2$ such that $p\in \ell\subset\Pi$, choose 
a point ${\Check p}$ and a line ${\Check \ell}$ in ${\Check \OO\PP}^2$ such that 
${\Check p}\in {\Check \ell} \subset{\Check \Pi}$, and finally a point $q\in \Lambda$. 
\bigskip

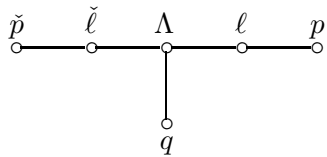
\begin{center}
\setlength{\unitlength}{5mm}
\begin{picture}(20,5)(-5,0)
\multiput(-.3,3.8)(2,0){5}{$\circ$}
\multiput(0,4)(2,0){4}{\line(1,0){1.7}}
\put(3.7,1.8){$\circ$}
\put(-.3,4.4){${\check p}$}
\put(1.7,4.4){${\check \ell}$}
\put(3.55,4.4){$\Lambda$}
\put(5.7,4.4){$\ell$}
\put(7.7,4.4){$p$}
\put(3.7,1.3){$q$}
\put(3.85,2.1){\line(0,1){1.8}}
\end{picture}
\end{center}

We call this data a complete $E_6$-flag. By \cite{tits}, there is a bijective 
correspondance between the set of Borel subgroups of $E_6$ and the set of complete 
$E_6$-flags: this is a direct generalization of the usual fact that a Borel subgroup of 
$SL_n$ is the stabilizer of a unique flag of vector subspaces of ${\bf C}^n$. 

We will not need this, but to complete the picture let us mention that 
the correspondance between $\Pi$, ${\Check \Pi}$ and $\Lambda$ can 
be described as follows:
$$\Pi = \bigcap_{z\in {\Check \Pi}}(T_z{\Check \OO}\PP^2)^{\perp} = 
\bigcap_{y\in \Lambda}y\JO.$$

\medskip This description of Borel subgroups will 
be useful to construct Schubert varieties in $\OO\PP^2$. Indeed, 
any subvariety of the Cayley plane that can defined in terms of a complete (or incomplete)
$E_6$-flag, must be a finite union of Schubert varieties. 

\smallskip Let us apply this principle in small codimension. The data ${\Check p},
{\Check \ell}, \Lambda$ 
from our $E_6$-flag are respectively a point, a line and a plane in ${\Check \OO\PP}^2$. 
They define special linear sections of $\OO\PP^2$, of respective codimensions $1$, $2$ and $3$. 
We read on the Hasse diagram that these sections are irreducible Schubert varieties. 

Something more interesting happens in codimension four, since we can read on the Hasse
diagram that a well-chosen codimension four linear section of $\OO\PP^2$ should split 
into the union of two Schubert varieties, of degrees $33$ and $45$. The most degenerate 
codimension four sections must correspond to very special $\PP^3$'s in ${\Check \OO\PP}^2$. 
We know from \cite{LMhom} that ${\Check \OO\PP}^2$ contains a whole family of $\PP^3$'s, in fact 
a homogeneous family parametrized by $E_6/P_{2,4}$. In terms of our $E_6$-flag, that 
means that a unique member of this family is defined by the pair $(q,\ell)$. 
 
We can describe explicitely a $\PP^3$ in $\OO\PP^2$ in the following way. Choose a 
non-zero vector $z\in \OO$, of zero norm. Then the closure of the set 
$$\Bigg\{\begin{pmatrix} 0 & x & 0 \\ \bar{x} & 0 & 0 \\ 0 & 0 & 0 \end{pmatrix},
\quad x\in L(z)\Bigg\},$$
is a three dimensional projective space $\PP^3_z$ in $\OO\PP^2$. Let us  take the
orthogonal of this space with respect to the quadratic form $Q(X)=\trace (X^2)$, 
and cut it with $\OO\PP^2$. We obtain two codimension $4$ subvarieties $Z_1$ 
and $Z_2$, respectively the closures of the following affine cells $Z_1^0$ and $Z_2^0$:
\begin{eqnarray}
Z_1^0 & = & \Bigg\{\begin{pmatrix} 1 & x & y \\ \bar{x} & 0 & \bar{x}y \\ 
\bar{y} & \bar{y}x & \bar{y}y \end{pmatrix}, \quad x\in L(z), 
y\in\OO\Bigg\}, \\
Z_2^0 & = & \Bigg\{\begin{pmatrix} 0 & u & uv \\ \bar{u} & 1 & v \\ 
\bar{v}\bar{u} & \bar{v} & \bar{v}v \end{pmatrix}, \quad u\in L(z), v\in\OO\Bigg\}.
\end{eqnarray}
The sum of the degrees of these two varieties is equal to $78$. The corresponding 
cycles are linear combinations of Schubert cycles with non negative coefficients.
But in codimension $4$ we have only two such 
cycles, $\s'_4$ and $\s_4''$, of respective degrees $33$ and $45$. The only 
possibility is that the cycles $[Z_1]$ and $[Z_2]$ coincide, up to the order, with 
$\s_4'$ and $\s_4''$. 

To decide which is which, let us cut $Z_1$ with $H_1=\{c_1=0\}$.

\begin{lemm}
The hyperplane section $Y_1=Z_1\cap H_1$ has two components $Y_{1,1}$ and $Y_{1,2}$.
One of these two components, say $Y_{1,1}$, is the closure of 
$$Y_{1,1}^0=\Bigg\{\begin{pmatrix} 0 & 0 & t \\ 0 & 0 & s \\ 
\bar{t} & \bar{s} & 1\end{pmatrix}, \quad  q(s)=q(t)=0, \; \bar{s}t=0\Bigg\}.$$
It is a cone over the spinor variety $\SS_{10}\subset\PP^{15}$. 
\end{lemm}

Recall that the spinor variety $\SS_{10}$ is one of the two families of maximal 
isotropic subspaces of a smooth eight-dimensional quadric.
Its appearance is not surprising, since we have 
seen on the weighted Dynkin diagram of $\OO\PP^2=E_6/P_6$ that the semi-simple
part of $P_6$ is a copy of ${\rm Spin}_{10}$. At a given point of $p\in\OO\PP^2$, 
the stabilizer $P_6$ and its subgroup ${\rm Spin}_{10}$ act on the tangent space, 
which is isomorphic as a ${\rm Spin}_{10}$-module to a half-spin representation,
say $\Delta_+$. From \cite{LMhom}, we know that the family of lines through $p$, that 
are contained in $\OO\PP^2$, is isomorphic to the spinor variety $\SS_{10}$, since
it is the closed ${\rm Spin}_{10}$-orbit in $\PP\Delta_+$.

In particular, to each point $p$ of $\OO\PP^2$ we can associate a subvariety, the union 
of lines through that point, which is a cone ${\mathcal C}(\SS_{10})$ over the 
spinor variety. This is precisely what is $Y_{1,1}$. Note that we get a Schubert variety
in the Cayley plane. Moreover, since we can choose a Borel subgroup of ${\rm Spin}_{10}$ 
inside a Borel subgroup of $E_6$ contained in $P_6$, we obtain a whole series of 
Schubert varieties which are isomorphic to cones over the Schubert subvarieties of 
$\SS_{10}$. These Schubert varieties can be described in terms of incidence relations
with an isotropic reference flag which in principle can be deduced from our reference
$E_6$-flag. 

The Hasse diagram of Schubert varieties in $\SS_{10}$ is the following:

\begin{center}
\setlength{\unitlength}{3.5mm}
\begin{picture}(30,12)(-8,2)

\multiput(-4.2,7.7)(2,0){3}{$\bullet$}
\multiput(1.8,9.7)(4,0){3}{$\bullet$}
\multiput(1.8,5.7)(4,0){3}{$\bullet$}
\multiput(3.8,7.7)(4,0){3}{$\bullet$}
\multiput(13.8,7.7)(2,0){2}{$\bullet$}
\put(7.8,11.7){$\bullet$}
\put(3.8,3.7){$\bullet$}

\multiput(-4.2,6.8)(2,0){3}{$1$}
\put(1.7,4.8){$1$}
\put(1.8,10.4){$1$}
\put(3.8,8.5){$2$}
\put(3.8,2.8){$1$}
\put(5.8,4.8){$3$}
\put(5.8,10.5){$2$}
\put(7.8,6.8){$5$}
\put(7.8,12.5){$2$}
\put(9.8,4.8){$5$}
\put(9.8,10.5){$7$}
\multiput(11.7,6.8)(2,0){3}{$12$}

\put(-4,8){\line(1,0){4}}
\put(0,8){\line(1,1){2}}
\put(0,8){\line(1,-1){4}}
\put(2,10){\line(1,-1){4}}
\put(6,10){\line(1,-1){4}}
\put(8,12){\line(1,-1){4}}
\put(2,6){\line(1,1){6}}
\put(4,4){\line(1,1){6}}
\put(10,6){\line(1,1){2}}
\put(12,8){\line(1,0){4}}

\end{picture}
\end{center}

The identification of the cone of lines in $\OO\PP^2$ through some given point, with 
the spinor variety $\SS_{10}$, is not so obvious. Consider the map 
$$\nu_2 : \OO\op\OO\ra\JZ, \qquad \nu_2(x,y)=\begin{pmatrix} x\bar{x}  &  x\bar{y} \\
y\bar{x}  &  y\bar{y}   \end{pmatrix}.$$
We want to identify $\PP\nu_2^{-1}(0)$ with $\SS_{10}$. The following result is due to P.E. Chaput:

\begin{prop}
Let $(x,y)\in\nu_2^{-1}(0)$. The image of the tangent map to $\nu_2$ at $(x,y)$ 
is a $5$-dimensional subspace of $\JZ$, which is isotropic with respect to the 
determinantal quadratic form on $\JZ$. Moreovoer, this induces an isomorphism 
between $\PP\nu_2^{-1}(0)$ and the spinor variety $\SS_{10}$.
\end{prop}

In fact, we can obtain this way the two families of maximal isotropic subspaces
in $\JZ$, just by switching the two diagonal coefficients in the definition of $\nu_2$. 
The spin group ${\rm Spin}_{10}$ can also be described very nicely.

\medskip
But let's come back to the Schubert varieties in $\SS_{10}$. Taking cones over 
them, we get Schubert subvarieties which define a subdiagram of the 
Hasse diagram of $\OO\PP^2$. We drew this subdiagram in thicklines on the picture below. \bigskip

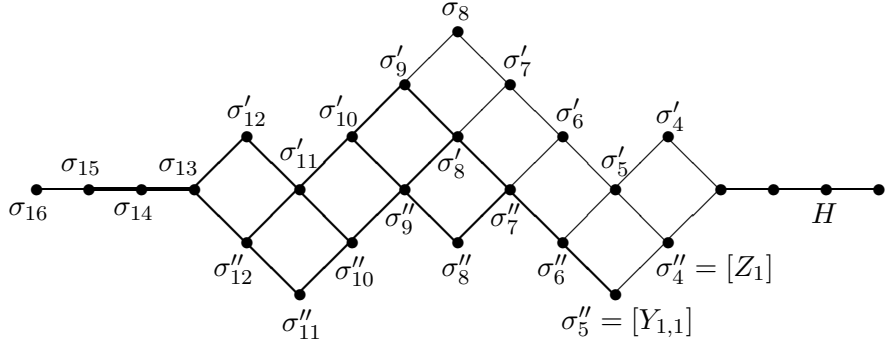
\begin{center}
\setlength{\unitlength}{3.5mm}
\begin{picture}(40,12)(-8,2)

\multiput(-6.2,7.7)(2,0){4}{$\bullet$}
\multiput(1.8,9.7)(4,0){5}{$\bullet$}
\multiput(1.8,5.7)(4,0){5}{$\bullet$}
\multiput(3.8,7.7)(4,0){5}{$\bullet$}
\multiput(3.8,3.7)(12,0){2}{$\bullet$}
\multiput(7.8,11.7)(4,0){2}{$\bullet$}
\multiput(21.8,7.7)(2,0){3}{$\bullet$}
\put(9.8,13.7){$\bullet$}
\put(.8,4.7){$\sigma_{12}''$}
\put(1.3,10.7){$\sigma_{12}'$}
\put(-1.3,8.7){$\sigma_{13}$}
\put(-3,7){$\sigma_{14}$}
\put(-5,8.7){$\sigma_{15}$}
\put(-7,7){$\sigma_{16}$}
\put(9.45,4.7){$\s_8''$}
\put(9.3,8.7){$\s_8'$}
\put(9.45,14.6){$\s_8$}
\put(23.5,6.7){$H$}
\put(17.6,4.7){$\s_4''=[Z_1]$}
\put(17.6,10.5){$\s_4'$}
\put(14,2.7){$\s_5''=[Y_{1,1}]$}
\put(15.5,8.8){$\s_5'$}
\put(13.1,4.7){$\s_6''$}
\put(13.8,10.6){$\s_6'$}
\put(11.3,6.6){$\s_7''$}
\put(11.8,12.6){$\s_7'$}
\put(7.3,6.6){$\s_9''$}
\put(7.1,12.6){$\s_9'$}
\put(5.4,4.6){$\s_{10}''$}
\put(4.8,10.7){$\s_{10}'$}
\put(3.4,2.6){$\s_{11}''$}
\put(3.3,9.1){$\s_{11}'$}

\put(8,12){\line(1,1){2}}
\put(10,10){\line(1,1){2}}

\put(10,14){\line(1,-1){8}}
\put(20,8){\line(1,0){6}}
\put(-6,8){\line(1,0){2}}
\put(18,10){\line(1,-1){2}}
\put(14,6){\line(1,1){4}}
\put(16,4){\line(1,1){4}}
\put(12,8){\line(1,1){2}}

\thicklines
\put(-4,8){\line(1,0){4}}
\put(0,8){\line(1,1){2}}
\put(0,8){\line(1,-1){4}}
\put(2,10){\line(1,-1){4}}
\put(6,10){\line(1,-1){4}}
\put(8,12){\line(1,-1){8}}
\put(2,6){\line(1,1){6}}
\put(4,4){\line(1,1){6}}
\put(10,6){\line(1,1){2}}

\end{picture}
\end{center}

We also indicate on this picture the indexing of Schubert classes that
we use in the sequel, rather than the indexing by the Weyl group. 

In principle, we are able to describe any of these Schubert varieties geometrically 
in terms of our reference $E_6$-flag.

\medskip\noindent {\it Proof of Lemma 4.1}. First note that $Y_1$ does not meet the 
two affine cells $\OO\PP^2_1$ and $\OO\PP^2_2$ (see section 2). Moreover, it is easy
to check that $Y_1\cap \OO\PP^2_{\infty}$ has dimension at most ten, hence strictly
smaller dimension than $Y_1$. Therefore, $Y_1$ is the closure of its intersection 
with $\OO\PP^2_3$, namely
$$Y_1\cap \OO\PP^2_3= 
\Bigg\{\begin{pmatrix} 0 & \bar{s}t & t \\ \bar{t}s & 0 & s \\ 
\bar{t} & \bar{s} & 1 \end{pmatrix}, \quad q(s)=q(t)=0, \; \bar{s}t\in L(z)\Bigg\}.$$
For a given non zero $s$, the product $\bar{s}t$ must belong to $L(z)\cap L(\bar{s})$,
the intersection of two maximal isotropic spaces of the same family. In particular, 
this intersection has even dimension. 

Generically, the intersection $L(z)\cap L(\bar{s})=0$, and we obtain
$$Y_{1,1}^0=\Bigg\{\begin{pmatrix} 0 & 0 & t \\ 0 & 0 & s \\ 
\bar{t} & \bar{s} & 1\end{pmatrix}, \quad  q(s)=q(t)=0, \; \bar{s}t=0\Bigg\}\subset Y_1.$$
We have seven parameters for $s$, and for each $s\ne 0$, $t$ must belong to $L(s)$, 
which gives four parameters. In particular, $Y_{1,1}^0$ is irreducible of dimension 
$11$, and its closure is an irreducible component of $Y_1$. 

The intersection $L(z)\cap L(\bar{s})$ has dimension two exactly when the line joining
$z$ to $\bar{s}$ is isotropic, which means that $\bar{s}$ belongs to the intersection of 
the quadric $q=0$ with its tangent hyperplane at $z$. This gives six parameters for $s$, 
and for each $s$, five parameters for $t$, which must be contained in the intersection of
the quadric with a six-dimensional linear space. Therefore, the closure of 
$$Y_{1,2}^0=\Bigg\{\begin{pmatrix} 0 & \bar{s}t & t \\ \bar{t}s & 0 & s \\ 
\bar{t} & \bar{s} & 1\end{pmatrix}, \quad  q(s)=q(t)=0, \; \dim L(z)\cap L(\bar{s})=2
\Bigg\},$$
is another component of $Y_1$. 

The remaining possibility is that $s$ be a multiple of $z$, but the corresponding subset of 
has dimension smaller than eleven. Hence $Y_1=Y_{1,1}\cup Y_{1,2}$. \qed

\medskip We conclude that $Z_1$ has degree $45$, while $Z_2$ has degree $33$. Indeed, 
if $Z_1$ had degree $33$, we would read on the Hasse diagram that its proper hyperplane 
sections are always irreducible, and we have just verified that this is not the case. 

Note that $Z_1$ and $Z_2$ look very similar at first sight. Nevertheless, a computation
similar  to the one we have just done shows that if we cut $Z_2$ by the hyperplane $H_2=
\{c_2=0\}$, we get an irreducible variety, the difference with $Z_1$ coming from the fact
that we now have to deal with maximal isotropic subspaces which are not on the same 
family. The difference between $Z_1$ and $Z_2$ is therefore just a question of spin...

\section{Intersection numbers}

We now determine the multiplicative structure of the Chow ring $A^*(\OO\PP^2)$. 
A priori, we have several interesting informations on that ring structure. We have
already seen in section 2 that the Pieri formula determines combinatorially the 
product with the hyperplane class. Another important property is that Poincar\'e 
duality has a very simple form in terms of Schubert cycles: the basis $(\s_w)_{w\in W^0}$ 
is, up to order, self-dual; more precisely its dual basis is $(\s_{w^*})_{w\in W^0}$, 
where the involution $w\mapsto w^*$ is very simple to define on the Hasse diagram:
it is just the symmetry with respect to the vertical line passing through the 
cycles of middle dimension. Finally, we know from Poincar\'e duality and general 
transversality arguments that any {\it effective} cycle must  be a linear combination of 
Schubert  cycles with non negative coefficients. 

\smallskip
This is the information we have on any rational homogeneous space. For what concerns 
the Cayley plane, we begin with an obvious observation:

\begin{prop} The Chow ring $A^*(\OO\PP^2)$ is generated by the hyperplane class
$H$, the class $\s_4'$, and the class $\s_8$ of an $\OO$-line. 
\end{prop}

More precisely, one can directly read on the Hasse diagram and from the Pieri 
formula that as a vector space, the Chow ring is generated by classes of type
$H^i$, $\s_4'H^j$ and $\s_8H^k$. For example, we have the  relations 
\begin{eqnarray}
H^4 & = & \s_4'+\s_4'', \\
\s_4' H^4 & = & \s_8+3\s_8'+2\s_8'', \\
\s_4'' H^4 & = & \s_8+4\s_8'+3\s_8'', \\
\s_8 H^4 & = & \s_{12}'+\sigma_{12}'', \\
\s_8' H^4 & = & 3\s_{12}'+4\sigma_{12}'', \\
\s_8'' H^4 & = & 2\s_{12}'+3\sigma_{12}''.
\end{eqnarray}
As a consequence, the multiplicative structure of the Chow ring will be completely 
determined once we'll have computed the intersection products $(\s^8)^2, 
\s_4'\s_8$
and $(\s_4')^2$. (Note that the Hasse diagram and the Pieri formula can be used to 
derive relations in dimension $9$ and $13$, but these relations are not sufficient to
determine the whole ring structure.)

\begin{prop} We have the following relations in the Chow ring:
\begin{eqnarray} \s_8^2 & = & 1, \\ \s_4' \s_8 & = & \s_{12}', \\ 
\s_4'' \s_8 & = & \s_{12}''.
\end{eqnarray}
\end{prop}

\proof Recall that $\s_8$ is the class of an $\OO$-line in $\OO\PP^2$, and that 
we know that the geometry of these lines is similar to the usual line geometry 
in $\PP^2$: namely, two generic lines meet transversely in one point. This implies
immediately that $\s_8^2=1$. 

To compute $\s_4' \s_8$ and $\s_4'' \s_8$, we cut the Schubert varieties $Z_1$ and 
$Z_2$ introduced in section 4, 
whose class we know to be $\s_4'$ and $\s_4''$, with the $\OO$-line $L$ defined 
in $\OO\PP^2$ by the conditions $x_1=x_2=r_3=0$. We get transverse intersections
$$\begin{array}{rcl}
Z_1\cap L & = & \Bigg\{\begin{pmatrix} r & y & 0 \\ \bar{y} & 0 & 0 \\ 
0 & 0 & 0\end{pmatrix}, \quad y\in L(z)\Bigg\}, \\ & & \\
Z_2\cap L & = & \Bigg\{\begin{pmatrix} 0 & y & 0 \\ \bar{y} & r & 0 \\ 
0 & 0 & 0\end{pmatrix}, \quad y\in L(z)\Bigg\}. 
\end{array}$$
These are two four dimensional projective spaces $\PP^4_1$ and $\PP^4_2$ inside 
$\OO\PP^2$, which look very similar. But there is actually a big difference: $\PP^4_1$ 
is extendable, but $\PP^4_2$ is not !
Indeed, a $\PP^5$ in $\OO\PP^2$ containing $\PP^4_1$ or $\PP^4_2$  must be of the form, 
respectively:
$$\Bigg\{\begin{pmatrix} r & y & s \\ \bar{y} & 0 & 0 \\ 
\bar{s} & 0 & 0\end{pmatrix}, \; y\in L(z)\Bigg\}, \quad {\rm and}\quad 
\Bigg\{\begin{pmatrix} 0 & y & 0 \\ \bar{y} & r & s \\ 
0 & \bar{s}  & 0\end{pmatrix}, \; y\in L(z)\Bigg\},$$
where $s$ describes some line in $\OO$. In the second case, the equation 
$sy=0$ must be verified identically, and we can take $s$ on the line 
${\bf C}{\bar z}$:
thus $\PP^4_2$ is extendable. But in the first case, we need the identity $y\bar{s}=0$
for all $y\in L(z)$, which would imply that $L(z)\subset R(s)$: this is impossible, 
and $\PP^4_1$ is not extendable. The proposition follows -- see the third observation at
the end of section 3. \qed

\medskip We now have enough information to complete the multiplication table. 
First, we know by Poincar\'e duality that 
\begin{eqnarray}
 (\s_8)^2  =  (\s_8')^2  =  (\s_8'')^2  =  1, \\
 \s_8\s_8'  =   \s_8'\s_8''   =  \s_8\s_8''  =  0, \\
  \s_4'\s_{12}'  = \s_4''\s_{12}''  =  1, \\
   \s_4'\s_{12}''  =  \s_4''\s_{12}'  =  0. 
\end{eqnarray}
Suppose that we have 
\begin{eqnarray}
\nonumber (\s_4')^2 & = & \mu_0\s_8+\mu_1\s_8'+\mu_2\s_8'', \\
\nonumber (\s_4'')^2& = &\nu_0\s_8+\nu_1\s_8'+\nu_2\s_8'', \\
\nonumber \s_4'\s_4'' & =& \g_0\s_8+\g_1\s_8'+\g_2\s_8'',
\end{eqnarray}
\noindent for some coefficients to be determined. Cutting with $\s_8$, we get $\mu_0=\nu_0=1$. 
The equations (3), (4), (5) give the relations 
$$\begin{array}{ccc}
\mu_0+\g_0=1, & \mu_1+\g_1=3, & \mu_2+\g_2=2, \\
\nu_0+\g_0=1, & \nu_1+\g_1=4, & \nu_2+\g_2=3.
\end{array}$$
In particular, $\g_0=0$. Now, we compute $(\s_4')^2(\s_4'')^2$ in two ways 
to obtain the relation
$$\g_1^2+\g_2^2=\mu_0\nu_0+\mu_1\nu_1+\mu_2\nu_2.$$
Eliminating the $\mu_i$'s and $\nu_i$'s, we get that $7\g_1+5\g_2=19$. But $\g_1$ and $\g_2$
are non negative integers, so the only possibility is that $\g_1=2$, $\g_2=1$. Thus:
\begin{eqnarray}
(\s_4')^2 & = & \s_8+\s_8'+\s_8'', \\
(\s_4'')^2 & = & \s_8+2\s_8'+2\s_8'', \\
\s_4'\s_4'' & = & 2\s_8'+\s_8''.
\end{eqnarray}
And this easily implies that 
\begin{eqnarray}
\s_4'\s_8' & = & \s_{12}'+2\s_{12}'', \\
\s_4'\s_8'' & = & \s_{12}'+\s_{12}'', \\
\s_4''\s_8' & = & 2\s_{12}'+2\s_{12}'', \\
\s_4''\s_8'' & = & \s_{12}'+2\s_{12}''.
\end{eqnarray}
\bigskip

\section{The Borel presentation}

We now turn to the Borel presentation of the Chow ring of $\OO\PP^2$. This is the 
ring isomorphism
$$A^*(\OO\PP^2)_{\QQ}\simeq \QQ [\cP]^{W_0}/\QQ [\cP]^W_+, $$
where $\QQ [\cP]^{W_0}$ denotes the ring of $W_0$-invariants polynomials on the 
weight lattice, and $\QQ [\cP]^W_+$ is the ideal of $\QQ [\cP]^{W_0}$ generated 
by $W$-invariants without constant term (see \cite{hiller}, ??). 

The ring $\QQ [\cP]^{W_0}$ is easily determined: it is generated by $\om_6$, and 
the subring of $W_0$ invariants in the weight lattice of ${\rm Spin}_{10}$. It is 
therefore the polynomial ring in the elementary symmetric functions 
$e_{2i}=c_i(\e_1,\ldots ,\e_5)$, $1\le i\le 4$, and in $e_5=\e_1\cdots\e_5$.

The invariants of $W$, the full Weyl group of $E_6$, are more difficult to determine, 
although we know their fundamental degrees. But since we know how to compute the 
intersection products of any two Schubert cycles, we just need to express the 
$W_0$-invariants in terms of the Schubert classes. This can be achieved, following
\cite{bgg}, by applying suitable difference operators to these invariants. 

Since we give a prominent role to the subsystem of $E_6$ of type $D_5$, it is natural to choose 
for the first five simple roots the usual simple roots of $D_5$, that is, in a euclidian 
$6$-dimensional space with orthonormal basis $\e_1,\ldots ,\e_6$, 
$$\begin{array}{rcl}
\a_1 & = & \e_1-\e_2, \\
\a_2 & = & \e_2-\e_3, \\
\a_3 & = & \e_3-\e_4, \\
\a_4 & = & \e_4-\e_5, \\
\a_5 & = & \e_4+\e_5, \\
\a_6 & = & -\frac{1}{2}(\e_1+\e_2+\e_3+\e_4+\e_5)+\frac{\sqrt{3}}{2}\e_6.
\end{array}$$
The fundamental weights are given by the dual basis:
\begin{eqnarray}
\nonumber\om_1 & = & \e_1+\frac{1}{\sqrt{3}}\e_6, \\
\nonumber\om_2 & = & \e_1+\e_2+\frac{2}{\sqrt{3}}\e_6, \\
\nonumber\om_3 & = & \e_1+\e_2+\e_3+\frac{3}{\sqrt{3}}\e_6, \\
\nonumber\om_4 & = & \frac{1}{2}(\e_1+\e_2+\e_3+\e_4-\e_5)+\frac{\sqrt{3}}{2}\e_6, \\
\nonumber\om_5 & = & \frac{1}{2}(\e_1+\e_2+\e_3+\e_4+\e_5)+\frac{5}{\sqrt{3}}\e_6, \\
\nonumber\om_6 & = & -\frac{1}{2}(\e_1+\e_2+\e_3+\e_4+\e_5)+\frac{\sqrt{3}}{2}\e_6.
\end{eqnarray}

The action of the fundamental reflexions on the weight lattice is specially
simple in the basis $\e_1,\ldots,\e_5,\a_6$. Indeed, $s_1,s_2,s_3$ and $s_4$ are
just the transpositions $(12), (23), (34), (45)$. The reflection $s_5$ affects
$\e_4,\e_5$ and $\a_6$, which are changed into $-\e_5,-\e_4$ and $\a_6+\e_4+\e_5$.
Finally, $s_6$ changes each $\e_i$ into $\e_i+\a_6/2$, and of course $\a_6$ into $-\a_6$. 

It is then reasonably simple to compute the corresponding divided differences with {\sc Maple}. 
We obtain: \vfill   \pagebreak

\begin{prop}
The fundamental $W_0$-invariants are given, in the Chow ring of the Cayley plane, 
in terms of Schubert cycles by;
\begin{eqnarray}
e_2 & = & -\frac{3}{4}H^2, \\
e_4 & = & -\frac{27}{8}\s_4'+\frac{21}{8}\s_4'', \\
e_5 & = & \frac{3}{16}\s_5'-\frac{21}{32}\s_5'', \\
e_6 & = & -\frac{27}{16}\s_6'+\frac{87}{32}\s_6'', \\
e_8 & = & \frac{21}{128}\s_8+\frac{291}{256}\s_8'-\frac{519}{256}\s_8''.
\end{eqnarray}
\end{prop}

\medskip
This allows to compute any product in the Borel presentation of the Chow ring
of $\OO\PP^2$. 
\medskip

\section{Chern classes of the normal bundle}

Let $\cN$ denote the normal bundle to the Cayley plane $\OO\PP^2\subset\PP\JO$. 
We want to compute its Chern classes. 

First note that the restriction of $\JO$ to the Levi part 
$L\simeq {\rm Spin}_{10}\times {\bf C}^*$ of the parabolic 
subgroup $P_6$ of $E_6$, is 
$$\JO_{|L}\simeq W_{\om_6}\op W_{\om_5-\om_6}\op W_{\om_1-\om_6}.$$
Indeed, there is certainly the line generated by the highest weight vector, 
which gives a stable line on which $L$ acts through the character $\om_6$. 
After $\om_6$, there is a in $\JO$ a unique highest weight, $\om_5-\om_6$, 
which generates a $16$-dimensional half-spin module. Finally, the lowest weight
of $\JO$ is $-\om_1$, whose highest $W_0$-conjugate is $\om_1-\om_6$ and 
generates a copy of the natural $10$-dimensional representation of ${\rm Spin}_{10}$. 
Since these three modules give $1+16+10=27$ dimensions, we have the full 
decomposition. 

Geometrically, this decomposition of $\JO$ must be interpreted as follows. We 
have chosen a point $p$ of $\OO\PP^2$, corresponding to the line $\hat{p}=W_{\om_6}$. 
The tangent space to $\OO\PP^2$ at that point is given by the factor $W_{\om_5-\om_6}$. 
(More precisely, only the affine tangent space $\Hat{T}$ is a well-defined $P_6$-submodule
of $\JO$, and it coincides with $W_{\om_6}\op W_{\om_5-\om_6}$.) The remaining 
term $W_{\om_1-\om_6}$ corresponds to the normal bundle. To be precise, if $\cN_p$ 
denotes the normal space to $\OO\PP^2$ at $p$, there is a 
canonical identification 
$$\cN_p\simeq Hom(\hat{p},\JO/\Hat{T})=Hom(W_{\om_6},W_{\om_1-\om_6}).$$
In other words, the normal bundle $\cN$ to $\OO\PP^2$ is the homogeneous bundle 
$\cE_{\om_1-2\om_6}$ defined by the irreducible $P_6$-module $W_{\om_1-2\om_6}$. 

Since $\om_1=\e_1+\frac{1}{2}\om_6$, the weights of the normal bundle are
the $\pm\e_i-\frac{3}{2}\om_6$, and its Chern class is
\begin{eqnarray}
\nonumber  c(\cN)  = & \prod_{i=1}^5(1+\e_i-\frac{3}{2}\om_6)(1-\e_i-\frac{3}{2}\om_6) \\
\nonumber    = & \prod_{i=1}^5\Big( (1+\frac{3}{2}H)^2-\e_i^2\Big) \\
\nonumber    = & \sum_{i=0}^5(-1)^i(1+\frac{3}{2}H)^{10-2i}e_{2i},
\end{eqnarray}
where $e_{10}=e_5^2$. We know how to express this in terms of Schubert classes, 
and the result is as follows.

\begin{prop}
In terms of Schubert cycles, the Chern classes of the normal bundle to 
$\OO\PP^2\subset\PP\JO$ are:
\begin{eqnarray}
\nonumber c_{1}(\cN) & = & 15H \\
\nonumber c_{2}(\cN) & = & 102H^2 \\
\nonumber c_{3}(\cN) & = & 414H^3\\
\nonumber c_{4}(\cN) & = & 1107\s_4'+1113\s_4'' \\
\nonumber c_{5}(\cN) & = & 2025\s_4' H+2079\s_4'' H \\
\nonumber c_{6}(\cN) & = & 5292\s_6'+8034\s_6''\\
\nonumber c_{7}(\cN) & = & 4698\s_6' H+7218\s_6'' H\\
\nonumber c_{8}(\cN) & = & 2751\s_8+9786\s_8'+7032\s_8''\\
\nonumber c_{9}(\cN) & = & 963\s_8H+3438\s_8'H+2466\s_8''H\\
\nonumber c_{10}(\cN) & = & 153\s_8H^2+549\s_8'H^2+387\s_8''H^2
\end{eqnarray}
\end{prop}

Note that as expected, we get integer coefficients, while the fundamental $W^0$-invariants
are only rational combinations of the Schubert cycles. This is a strong indication that our
computations are correct.

\section{The final computation}

We want to compute the degree of the variety of reductions $Y_8$ introduced 
in \cite{IM}: we refer to this paper for the definitions, notations, and the 
proof of the facts we use in this section. This variety $Y_8$ is a smooth projective 
variety of dimension $24$, embedded in $\PP^{272}$. A $\PP^1$-bundle $Z_8$ over $Y_8$ 
can be identified with the blow-up of the projected Cayley plane $\overline{X_8}$ in 
$\PP\JO_0$, the projective space of trace zero Hermitian matrices of order three,
with coefficients in the Cayley octonions. 

Let $H$ denote the pull-back to $Z_8$ of the hyperplane class of $\PP\JO_0$, and $E$
the exceptional divisor of the blow-up. We want to compute 
$${\rm deg}\, Y_8=H(3H-E)^{24}.$$
We use the fact that the Chow ring of the exceptional divisor $E\subset Z_8$, since it is
the projectivization of the normal  is the quotient of the ring  $A^*(\OO\PP^2)[e]$ 
by the relation given by the Chern classes of the normal bundle $\Bar{\cN}$ of $\overline{X_8}$, 
namely
$$e^{9}+\sum_{i=1}^{9} (-1)^i c_i(\Bar{\cN})e^{9-i}=0.$$
The normal bundle $\Bar{\cN}$ of $\overline{X_8}$ is related to the normal bundle $\cN$ 
of $X_8=\OO\PP^2$ by an exact sequence $0\ra \cO(1)\ra\cN\ra
\Bar{\cN}\ra 0$, from which we can compute the Chern classes of $\Bar{\cN}$:

\begin{eqnarray}
\nonumber c_{1}({\Bar\cN}) & = & 14H \\
\nonumber c_{2}({\Bar\cN}) & = & 88H^2 \\
\nonumber c_{3}({\Bar\cN}) & = & 326H^3\\
\nonumber c_{4}({\Bar\cN}) & = & 781\s_4'+787\s_4'' \\
\nonumber c_{5}({\Bar\cN}) & = & 1244\s_4' H+1292\s_4'' H=2536\s_5'+1292\s_5'' \\
\nonumber c_{6}({\Bar\cN}) & = & 2756\s_6'+4206\s_6''\\
\nonumber c_{7}({\Bar\cN}) & = & 1942\s_6' H+3012\s_6'' H=1942\s_7'+4954\s_7''\\
\nonumber c_{8}({\Bar\cN}) & = & 809\s_8+2890\s_8'+2078\s_8''\\
\nonumber c_{9}({\Bar\cN}) & = & 154\s_8H+548\s_8'H+388\s_8''H=702\s_9'+936\s_9''\\
\nonumber c_{10}({\Bar\cN}) & = & -\s_8H^2+\s_8'H^2-\s_8''H^2=0 \;  !
\end{eqnarray}
The fact that we get $c_{10}({\Bar\cN})=0$, which must hold since ${\Bar\cN}$
has rank $9$, is again a strong indication that we did no mistake. 

\medskip
To complete our computation, we must compute the intersection products $H^{25-i}E^i$ in the 
Chow ring of $Z_8$. For $i>0$, this can be computed on the exceptional divisor; since the 
restriction of the class $E$ to the exceptional divisor is just the relative hyperplane section,
that is, the class $e$, we have $H^{25-i}E^i=H^{25-i}e^{i-1}$, the later product being computed 
in $A^*(E)$. We still denoted by $H$ the pull-back of the hyperplane section from $\OO\PP^2$. 

\begin{lemm} Let $\s\in A^{16-k}(\OO\PP^2)$. Then $\s e^{8+k}=\s s_k({\Bar\cN})$, where
$s_k({\Bar\cN})$ denotes the $k$-th Segre class of the normal bundle ${\Bar\cN}$. The former
product is computed in $A^*(E)$, and the the later in $A^*(\OO\PP^2)$.
\end{lemm}

\proof Induction, using the relation $e^{9}+\sum_{i=1}^{9} (-1)^ic_i(\Bar{\cN})e^{9-i}=0$, and the 
fact that the Segre classes are related to the Chern classes by the formally similar relation 
$s_k({\Bar\cN})+\sum_{i=1}^{9} (-1)^ic_i(\Bar{\cN})s_{k-i}({\Bar\cN})=0$.\qed

\medskip 
We use the later relation to determine the Segre classes inductively. We obtain
\begin{eqnarray}
\nonumber s_1({\Bar\cN})& = & 14H, \\
\nonumber s_2({\Bar\cN})& = & 108H^2, \\
\nonumber s_3({\Bar\cN})& = & 606H^3, \\
\nonumber s_4({\Bar\cN})& = & 2763\s_4'+2757\s_4'', \\
\nonumber s_5({\Bar\cN})& = & 21624\s_5'+10752\s_5''  \\
\nonumber s_6({\Bar\cN})& = & 75492\s_6'+112602\s_6'',\\
\nonumber s_7({\Bar\cN})& = & 240534\s_7'+596598\s_7'', \\
\nonumber s_8({\Bar\cN})& = & 711489\s_8+2462397\s_8'+1750947\s_8'', \\
\nonumber s_9({\Bar\cN})& = & 8768196\s_9'+11600304\s_9'', \\
\nonumber s_{10}({\Bar\cN})& = & 53127900\s_{10}'+30193704\s_{10}'', \\
\nonumber s_{11}({\Bar\cN})& = & 206857602\s_{11}'+74823228\s_{11}'', \\
\nonumber s_{12}({\Bar\cN})& = & 491985531\s_{12}'+669523221\s_{12}'', 
\end{eqnarray}

\pagebreak

\begin{eqnarray}
\nonumber s_{13}({\Bar\cN})& = & 2657712312\s_{13}, \\
\nonumber s_{14}({\Bar\cN})& = & 5875513812\s_{14} \\
\nonumber s_{15}({\Bar\cN})& = & 12591161406\s_{15}.
\end{eqnarray}

\bigskip
This immediately gives the degree of $Y_8$, 
$${\rm deg}\,Y_8 = 3^{24}+\sum_{k=9}^{24}(-1)^k\binom{24}{k}3^{24-k}H^{25-k}s_{k-9}({\Bar\cN}).$$

\medskip

\begin{theo}
The degree of the variety of reductions $Y_8$ is 
$${\rm deg}\,Y_8 = 1\,047\,361\,761.$$
\end{theo}

\bigskip\noindent {\it Acknowledgements}. We warmly thank D. Markushevitch for writing the
Macaulay script which allowed the computation of the Segre classes above, and P.E. Chaput 
for useful duscussions.

\bigskip\noindent 
Atanas ILIEV, Institut de Math\'ematiques, Acad\'emie des Sciences 
de Bulgarie, 8 Rue Acad. G. Bonchev, 1113 Sofia, Bulgarie.

\noindent Email : ailiev@math.bas.bg

\smallskip\noindent 
Laurent MANIVEL, Institut Fourier, Laboratoire de Math\'ematiques, 
UMR 5582 (UJF-CNRS), BP 74, 38402 St Martin d'H\`eres Cedex, France.

\noindent Email : Laurent.Manivel@ujf-grenoble.fr

\end{document}